\documentclass{article}
\usepackage{amssymb, amsmath, latexsym, amscd, amsthm}
\usepackage{graphicx}

\newtheorem{theorem}{Theorem}[section]
\newtheorem{corollary}[theorem]{Corollary}
\newtheorem{lemma}[theorem]{Lemma}

\newtheorem{proposition}[theorem]{Proposition}

\newtheorem{remark}[theorem]{Remark}

\newcommand{\Q}{{\mathbb Q}} 
\newcommand{\g}{\mathcal G}
\newcommand{\h}{\mathcal {HG}}
\newcommand{\cv}{ \mathcal C} 
\newcommand{\F}{\mathcal F}
\newcommand{\bdry}{\partial}

\title{Cut vertices in commutative graphs}
\author{James Conant\footnote{E-mail: jconant@math.utk.edu}\\
{\it(Department of Mathematics, University of Tennessee at Knoxville,}\\ \it{Knoxville, TN 37996-1300)}\\
\\ Ferenc Gerlits\footnote{E-mail: fgerlits@renyu.hu}\\
{\it(Alfr\'ed R\'enyi Institute of Mathematics, Hungarian Academy of Sciences,}\\
{\it Re\'altanoda utca 13-15 H-1053, Budapest, Hungary)}\\
\\ Karen Vogtmann\footnote{E-mail: vogtmann@math.cornell.edu}\\
{\it(Department of Mathematics, Cornell University,}\\
{\it Ithaca, NY 14853-4201)}}

\date{}

\begin{document}
\maketitle

\begin{abstract}
The homology of Kontsevich's commutative graph complex parameterizes finite type invariants of odd dimensional manifolds. This {\it graph homology} is also  the twisted homology of Outer Space modulo its boundary, so gives
a nice point of contact between geometric group theory and quantum topology.
In this paper we give two different proofs (one algebraic, one geometric) that the commutative graph complex is quasi-isomorphic to the quotient complex obtained by modding out by graphs with cut vertices.  This quotient complex has the advantage of being smaller and hence
more practical for computations.  In addition, it supports a Lie bialgebra structure coming from a bracket
and cobracket we defined in a previous paper. As an application, we compute
the rational homology groups of the commutative graph complex up to rank $7$.
\end{abstract}

\large

\section{Introduction}
Graph homology was introduced by Kontsevich \cite{Ko1, Ko2}, who showed  that it computes the homology of a certain infinite dimensional 
Lie algebra $c_\infty$, and also 
parameterizes invariants of certain odd dimensional manifolds. 
 The best understood of these invariants are those associated to (rational) homology three spheres.  These are known as ``finite type" invariants,  and are analogs of the Goussarov-Vassiliev knot invariants. Alternate constructions of these finite type invariants have been found by
Le, Murakami and Ohtsuki \cite{lmo}, Kuperberg and Thurston \cite{kt}, 
and Bar-Natan, Garoufalidis, Rozansky and Thurston \cite{bgrt}.

Graph homology has a very simple definition.
The degree $k$ term of the graph complex $\g$  is spanned (over a field of characteristic $0$)
 by connected, ``oriented" graphs with
$k$ vertices, and the boundary operator $\partial\colon
\g_k\to\g_{k-1}$ is defined on a graph $G$  by adding together 
 all oriented graphs which can be obtained from $G$ by collapsing a single edge.Ê The notion of orientation
is the most subtle part of the whole story (see \cite[Sec. 2.3.1]{exposition}  for the many
equivalent notions), but suffice it to say that
it guarantees that $\partial^2=0$.
Graph homology is then the homology of this complex.

Graph homology is a special case of a more
general construction, where the graph complex is spanned by graphs decorated at each vertex
by an element of a cyclic operad (see \cite{exposition}).
Ordinary graph homology corresponds to the commutative operad.Ê Two
other important examples, also studied by
Kontsevich \cite{Ko1,Ko2}, are obtained by using the associative operad and the Lie
operad; the homology of the resulting
graph complexes is closely related to the cohomology of mapping class
groups of punctured surfaces and the
cohomology of outer automorphisms of free groups, respectively. 

Each general graph complex $\g$Ê may be considered as the primitive part of a
graded Hopf algebra $\h$, where
the product on $\h$ is given by disjoint union.ÊÊ In \cite{bialgebra}
we introduced a Lie bracket and cobracket on $\h$. These do not form
a compatible bialgebra structure on $\h$, and they do not restrict to
$\g$.Ê However, in \cite{bialgebra}
we also introduced the subcomplex $\mathcal B$ of $\g$ spanned by
connected graphs with no separating edges, and showed
that the Lie bracket and cobracket do restrict to $\mathcal B$;
furthermore, they are compatible, so giveÊ $\mathcal B$ aÊ Lie
bialgebra
structure. In the associative and Lie cases, the subcomplex $\mathcal
B$ isÊ quasi-isomorphic to $\g$, but this is not true in the
commutative case:
$\mathcal B$ and $\g$ do not have the same homology.

\begin{figure}
\begin{center}
\includegraphics[height=3cm]{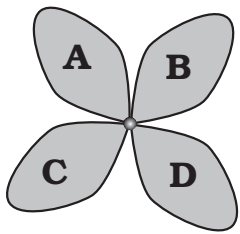}
\end{center}
\caption{Schematic for a cut vertex.} \label{cv1}
\end{figure}
In this paper we use a different approach in the commutative case to find a smaller chain complex
quasi-isomorphic to $\g$  which carries a Lie
bialgebra structure.Ê Specifically, we consider the subcomplex $\cv$ of
$\g$ spanned by graphs with at least one \emph{cut
vertex}, where a cut vertex is defined as a vertex whose deletion 
disconnects the graph (Figure~\ref{cv1}).
In Section 2 we use a spectral sequence argument to prove 

\begin{theorem}\label{quasiisometric}
The quotient map of chain complexes $ \g \to \g/\cv$ is an isomorphism
on homology.
\end{theorem}

In section 3 we recall the geometric interpretation of graph homology in terms of
Outer space from \cite{exposition}, and reprove Theorem~\ref{quasiisometric} from this point of view. 
Specifically, we show that the standard deformation retraction of Outer
space onto the subspace of graphs with no separating edges
 extends to the Bestvina-Feign bordification of Outer Space, and we
show that the image of the points at infinity consists precisely of the
closure of the space of graphs with cut vertices.
This deformation retraction induces a homology isomorphism on certain twisted
chain complexes, which can be identified with $\g$ and $\g/{\mathcal C}$.

In section 4 we recall the definition of the Lie bracket and cobracket, and
prove

\begin{theorem}\label{bialgebra} The Lie bracket and cobracket on $\h$
induce a compatible graded Lie bialgebra structure on $\g/\cv$.
\end{theorem}

{\bf Remark.} Despite the fact that the entire graph complex $\h$ does
not support a Lie bialgebra structure, Wee-Liang Gan \cite{gan} has recently shown that it
supports a strongly homotopy Lie bialegbra structure, and that this reduces to our
Lie bialgebra structure when one mods out by graphs with cut vertices.
\medskip

Finally, in the last section we exploit the fact that the quotient complex $\g/{\mathcal C}$ is
smaller than $\g$ to do some computer-aided calculations of graph homology.   Specifically, the elimination of cut vertices reduces the size of the vector spaces involved by about $30\%$, allowing us to calculate graph homology up to rank 7.

\section{Graphs with cut vertices}

Let $G$ be an orientedÊ graph with no separating edges. An oriented graph $H$ is said to {\it
retract to $G$} if
$G$ can be obtained from $H$ by collapsing each separating edge of $H$
to a point.
Denote by ${\mathcal R}^G$
the subspace of
$\g$ spanned by all graphs $H$ which retract to $G$.
Notice that according to our definitions, $\mathcal R^G$ is $1$ dimensional 
(with basis $G$) unless
$G$ has a cut vertex.
Ê Define a boundary operator
$\partial_s\colon{\mathcal R}^G_k\to{\mathcal R}^G_{k-1}$ by
$$\partial_s(G)=\sum_{e\,\,\, \text{separating}} G_e$$

\begin{figure}
\begin{center}
\includegraphics[height=7cm]{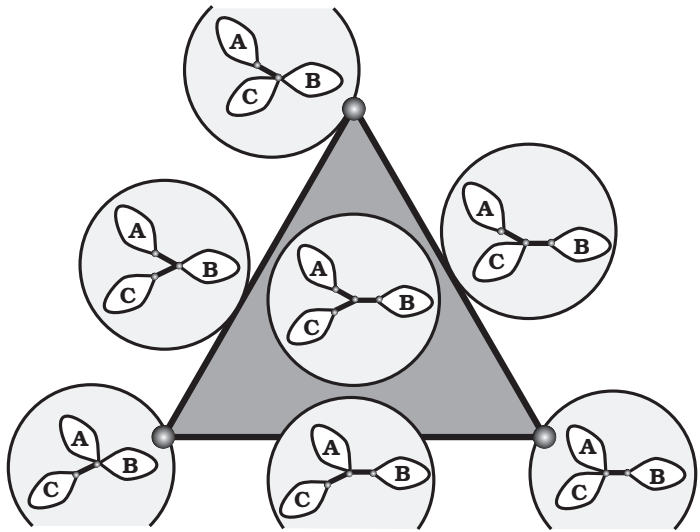} 
\end{center}
\caption{The simplicial realization of ${\mathcal R}^v_+$.}\label{cv0}
\end{figure}

\begin{lemma} \label{mainlem}
LetÊ $G$ be a connected graph with at least one cut vertex, but no
separating edges.Ê Then $({\mathcal R}^G,\partial_s)$ is an acyclic
complex.
\end{lemma}

\begin{proof}
First we prove the lemma under the assumption that $G$ has no automorphisms.
Fix a cut vertex $v$ of $G$, and let $c_1,\ldots,c_l$Ê be the connected
components
of $G\setminus\{v\}$.Ê Let ${\mathcal R}^{v}$ denote the subcomplex of
${\mathcal R}^G$ spanned by
graphs whose separating edges form a tree which collapses to $v$.
Grade $\mathcal R^v$ by the number of edges in the tree that collapses to $v$.
 We will show that the positive degree part of this subcomplex, ${\mathcal
R}^v_+$, is the simplicial homology of a contractible simplicial complex,
and that the map
$\partial_s\colon {\mathcal R^v_1\to\mathcal R^v_0}$ serves as an
augmentation map, implying the complex is acyclic. A low degree example
is given in  Figure \ref{cv0}, which shows that the simplicial
realization of
${\mathcal R^v_+}$ where $v$ is a cut vertex that cuts the graph into three 
components is in fact a $2$-simplex. Notice that face maps correspond exactly to 
collapsing edges. Finally one needs to consider orientation.
One notion of orientation for commutative graphs is an orientation
of the vector space $\mathbb R^{\{edges\}}\oplus H_1(Graph;\mathbb R)$.
The first factor corresponds exactly to an orientation of the above 2-simplex,
whereas the orientation of $H_1$ gets carried along for the ride.

In general the simplicial structure of the blow ups of a vertex is easiest to
understand by considering sets of \emph{compatible} partitions of $\{c_i\}$
as opposed to trees.
IfÊ $H$ is a graph in ${\mathcal R}^{v}$ andÊ $T$ is the tree in $H$
which collapses to $v$,
then each edge
$e$ of
$T$ partitions the (preimages in $H$ of the) components
$c_i$ into two disjoint sets.Ê Different edges correspond to different
partitions, which are {\it
compatible} in the following sense:Ê If
$P_1=X_1\cup Y_1$ and $P_2 = X_2\cup Y_2$, then either
$$X_1\subset X_2,\,\, X_1\subset Y_2,\,\, Y_1\subset X_2,\,\, Y_1\subset
Y_2.$$ Conversely, any set of pairwise compatible partitions determines a
pair
$(H,T)$ which collapses to $(G,v)$. Figure~\ref{cv2} shows a tree that blows
up a vertex and two compatible partitions corresponding to two edges of the tree.

\begin{figure}
\begin{center}
\includegraphics[height=4cm]{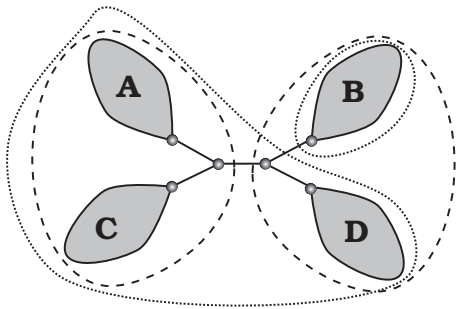}
\end{center}
\caption{Two compatible partitions of the set $\{A,B,C,D\}$
corresponding to two edges of a tree.}\label{cv2}
\end{figure}

ÊÊÊ Note that ${\mathcal R}^v_0$ is 1-dimensional, spanned by $G$.Ê Thus
Ê ${\mathcal R}^v$ is the augmented chain complex of a
simplicial complex whose vertices are partitions of the set
$\{c_1,\ldots,c_l\}$.Ê A set of $k+1$ partitions
forms a
$k$-simplex ifÊ partitions in the set are pairwise compatible.Ê The
partition $P$ separating $c_1$
from all other componentsÊ is compatible with every other set of partitions;
thus this simplicial complex isÊÊ a
cone on the vertex $P$, and is therefore contractible.Ê Thus the
augmented chain complex $({\mathcal R}^v,\partial_s)$
is acyclic.

If we grade the entire complex $({\mathcal R}^G,\partial_s)$ by the number of edges
in the forest formed by all separating edges, then it
 is the tensor product
of the complexes $({\mathcal R}^v,\partial_s)$ for
cut vertices $v$ of $G$. Thus ${\mathcal R}^G$ is acyclic.

Now let us return to the general case when the graph has a nontrivial automorphism group.
Let $\widehat{\mathcal R}^G$ be the chain complex of graphs obtained from ${\mathcal R}^G$ by
 distinguishing all edges and vertices in each graph.
(This kills off automorphisms.)
Note that $Aut(G)$ acts on $\widehat{\mathcal R}^G$, and that $\widehat{\mathcal R}^G/ Aut(G)\cong  {\mathcal R}^G$.
Over the rationals, finite groups have no homology, a fact which implies that
the chain complexes $\widehat{\mathcal R}^G$ and ${\mathcal R}^G$ are rationally quasi-isomorphic.
Now we are back in the case when graphs have no automorphisms, and we are done.
\end{proof}

\begin{theorem} The subcomplex $\cv$ of $\g$ spanned by graphs with at
least one cut vertex is acyclic.
\end{theorem}

\begin{proof}
FilterÊ $\cv$ by the number of separating edges in the graph;
i.e. let $\F_p\cv$ be the subcomplex of $\cv$ spanned by graphs with at
most $p$ separating edges.
Then
$$\F_0 \cv\subset \F_1\cv\subset \F_2\cv\ldots.$$

The boundary operator $\bdry\colon \cv_k\to\cv_{k-1}$ is the sum of two
boundary operators $\bdry_s$ and
$\bdry_{ns}$, where $\bdry_s$ collapses only separating edges, and
$\bdry_{ns}$ collapses only
non-separating edges.Ê These two boundary operators make $\cv$ into the
total complex of a double complex
$E^0$, where
$E^0_{p,q}=\F_p(\cv_{p+q})$, the vertical arrows are given by $\bdry_s$
and the horizontalÊÊ arrowsÊ by
$\bdry_{ns}$:
$$\begin{matrix}
&\downarrow&&\downarrow\\
\leftarrow &\F_p(\cv_{p+q-1})&\buildrel
\bdry_{ns}\over\leftarrow&\F_p(\cv_{p+q})&\leftarrow \\
&\ \ \downarrow_{\bdry_{s}}&&\ \ \downarrow_{\bdry_{s}}\\
\leftarrow &\F_{p-1}(\cv_{p+q-2})&\buildrel
\bdry_{ns}\over\leftarrow&\F_{p-1}(\cv_{p+q-1})&\leftarrow\\
&\downarrow&&\downarrow\\
Ê \end{matrix}
$$
Ê Note
that a graph in
$\F_p(\cv)$ has at least $p+1$ vertices, so that
$E^0_{p,q}=\F_p(\cv_{p+q})=0$ for $q<1$, and the double
complex is a first quadrant double complex.

We consider the spectral sequence associated to the vertical filtration
of this double complex.ÊÊ This
spectral sequence converges to the homology of the total complex $\cv$.
The
$E^1_{p,q}$ term is equal to $H_p(E^0_{*,q},\bdry_s)$, i.e. the $p$-th
homology of the $q$-th column.

For
each $q$, the column $E^0_{*,q}$ breaks up into a direct sum of chain
complexes $E^G_*$, one for each
graph
$G$ with $q$ vertices (at least one of which is a cut vertex) and no
separating edges.
A graph in $\F_p(\cv_{p+q})$ is in $E^G_*$ if
$G$ is the result of collapsing all of its separating edges, i.e.
$E^G_*={\mathcal R}^G$.
By Lemma \ref{mainlem}, ${\mathcal R}^G$ has no homology, so that
$E^1_{p,q}=0$ for all $p$ and $q$, and the complex $\cv$ is acyclic.
\end{proof}

Theorem \ref{quasiisometric} now follows immediately by the long exact
homology sequence of the pair $(\g,\cv)$.

\section{Geometric Interpretation, in terms of Outer Space}

In this section we will sketch a geometric proof of the main theorem.
This proof relies on the identification of the graph homology chain complex
with a twisted relative chain complex for Outer space, as described in \cite{exposition},
and also on a generalization of the {\it Borel-Serre bordification} of Outer space
defined by Bestvina and Feighn \cite{BF}.  A similar generalization is mentioned as a remark in their
paper, but details of proofs are not worked out.

Recall that Outer space $X_n$
is a  topological space which parameterizes finite metric graphs with
(free) fundamental group of rank $n$ (see \cite{survey}).ÊÊ Outer space
can be decomposed as a union of open simplices, and there are  several ways
to add a boundary to this space. The simplest is to formally add the union of
all missing faces
 to obtain a simplicial complex
$\overline{X}_n$, called the {\it simplicial closure of Outer Space}.
The bordification is more subtle; it is a
blown-up version of $\overline{X}_n$, which we will denote $\widehat{X}_n$.
The interiors of $\overline X_n$ and $\widehat X_n$ are both homeomorphic to
$X_n$, and the action of $Out(F_n)$ extends to the boundaries
$\partial
\overline{X}_n$ and $\partial
\widehat{X}_n$.  There is a natural quotient map
$q\colon\widehat X_n
\to 
\overline X_n$, which  is a homeomorphism on the interiors and in general has contractible point inverses.

In \cite{exposition} we showed
that the subcomplex 
${\mathcal G}^{(n)}$ of the graph complex $\g$ spanned by graphs with fundamental
group of  rank
$n$ can be identified with the relative chains on
$(\overline X_n,\bdry\overline X_n)$, twisted by the non-trivial {\it determinant}
action of $Out(F_n)$ on ${\mathbb R}$. Blowing up the boundary does not
change this picture;  
${\mathcal G}^{(n)}$ is also identified with the relative chains
on
$(\widehat{X}_n,\partial \widehat{X}_n)$, twisted by the same non-trivial 
action of $Out(F_n)$ on ${\mathbb R}$.

In this section we define an
equivariant deformation retraction $\widehat X_n \to \widehat Y_n$,
where $Y_n$ is the subspace of $X_n$ consisting of graphs with no
separating edges, and $\widehat Y_n$ denotes the closure of $Y_n$ in
$\widehat X_n$. The image  of $\partial \widehat{X}_n$ under this retraction, denoted $\widehat Z_n$,
 is the union of $\partial\widehat Y_n$ and 
the set of graphs with a cut vertex but no separating edges.Ê The deformation retraction
induces an isomorphism
$$C_*(\widehat X_n,\partial\widehat X_n)\otimes_{Out(F_n)}{\mathbb R}\cong
C_*(\widehat Y_n, \widehat Z_n)\otimes_{Out(F_n)}{\mathbb R}.$$
Tracing through the identification of $C_*(\widehat
X_n,\partial\widehat X_n)\otimes_{Out(F_n)}{\mathbb R}$ with ${\mathcal G}^{(n)},$ we see that the chains
$$C_*(\widehat Y_n, \widehat Z_n)\otimes_{Out(F_n)}{\mathbb R}$$ are identified with ${\mathcal
G}^{(n)}_{ns}/{\mathcal C}^{(n)}_{ns}$, where the subscript $ns$  denotes the subcomplex spanned by
graphs with no separating edges. 
 Since all graphs with separating edges also have cut vertices, this is
naturally isomorphic to
${\mathcal G}^{(n)}/{\mathcal C}^{(n)}$. This completes the sketch of the proof of Theorem~\ref{quasiisometric}, modulo the definition of
the bordfication and the retraction. The remainder of the section is devoted to just that.

\begin{figure}
\begin{center}
\includegraphics[height=4cm]{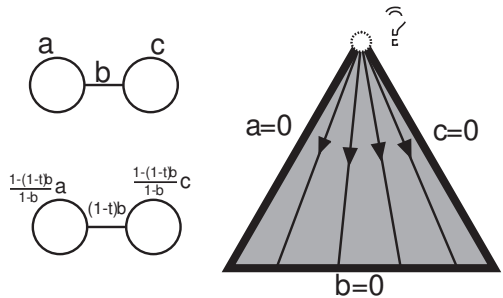}
\end{center}
\caption{Why the simplicial closure doesn't work.} \label{triangle0}
\end{figure}

It has long been known that $Y_n$ is an equivariant deformation retract
of $X_n$, but the deformation retraction, which uniformly shrinks
all separating edges while uniformly expanding all other edges, does not
extend to $\overline X_n$.Ê One can see this even for $n=2$ by
considering the 2-simplex corresponding to the ``barbell" graph (see
Figure \ref{triangle0}).Ê The deformation retraction sends each horizontal slice
linearly onto the bottom edge of the triangle, so that the deformation
cannot be extended continuously to the top vertex of the closed
triangle.

This difficulty can be resolved  by blowing
up the vertex of the triangle to a line, which records the (constant) ratio of the lengths of
the two loops of the ``barbell" graph along a geodesic in
$X_2$ coming into the vertex (see Fig. \ref{triangle0}).
This is the idea of 
the Borel-Serre bordification $\widehat X_n$.  Similar ideas are also
used in the Fulton-MacPherson and
Axelrod-Singer compactifications of configuration spaces of points in a manifold. 

To describe $\widehat X_n$ in general, we need the notion of 
a \emph{core graph}, which is defined to be 
a (not necessarily connected) graph with no separating edges and no
vertices of valence 0 or 1.  Every graph has a unique maximal core
subgraph,  called its
\emph{core}.

A point in $X_n$ is a marked, nondegenerate metric graph of rank $n$ with
total volume  $1$, where ``nondegenerate" means no edge is assigned the length $0$.
In the bordification $\widehat{X}_n$, we allow edges of a core subgraph to have length $0$, but in this
case there is a secondary metric, also of volume 1, given on the core subgraph.  The secondary metric may
also be zero on a smaller core subgraph, in which case there is a third metric of volume 1 on that core
subgraph, etc.  

 In general, a point of
$\widehat X_n$ consists of marked metric graph
$\Gamma_0$ and a properly nested (possibly empty) sequence 
$$\Gamma_1\supset \Gamma_2\supset\ldots\supset \Gamma_k$$
of core
subgraphs of $\Gamma_0$. Each $\Gamma_i$ is equipped with a metric of volume 1; the metric on $\Gamma_0$ is the {\it
primary metric}, the metric on $\Gamma_1$ the {\it secondary metric} etc.  Each
$\Gamma_i$ is the subgraph of
$\Gamma_{i-1}$ spanned by all edges of length $0$, and the chain is {\it nondegenerate} in the sense that
every edge of $\Gamma_0$ has non-zero length in exactly one $\Gamma_i$.   The space
$\widehat X_n$ is stratified as a union of open cells; the dimension of the cell containing 
$x=(\Gamma_0\supset\Gamma_1\supset\ldots\supset\Gamma_k)$ is $e(\Gamma_0)-k-1$, where $e(\Gamma_0)$ is the
number of edges of
$\Gamma_0$.  

\begin{figure}
\begin{center}
\includegraphics[height=8cm]{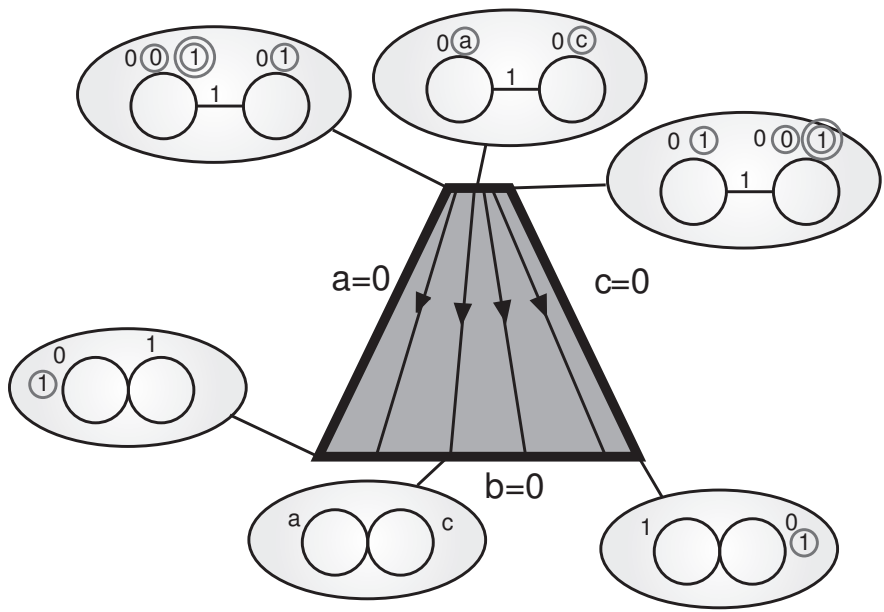}
\end{center}
\caption{The deformation retraction on a cell of the bordification.} \label{triangle}
\end{figure}

This construction is illustrated for $n=2$ in Fig. \ref{triangle}. In this figure, the number
of circles surrounding an edge length corresponds to the hierarchy of metrics.  A sequence of graphs on
which the volume of a core subgraph is shrinking to zero will approach a point on the boundary
which depends on the relative lengths of edges in the  core subgraph. If the metric is
shrinking uniformly on the core subgraph, then the limit is the graph whose primary metric vanishes on the
core subgraph, and where the secondary metric on the core subgraph is a rescaled version of the metric
restricted to the shrinking core subgraph. If parts of the core subgraph are shrinking at a faster rate than
others, the sequence will land in a face of higher codimension.

Bestvina and Feighn prove the following theorem.
 
\begin{theorem}[Bestvina-Feighn]\label{bfthm}
$\widehat{Y}_n$ is contractible, and the $Out(F_n)$ action on the interior extends to the whole space.  
\end{theorem}

The following theorem shows that $\widehat X_n$ is also contractible, and identifies the image of the boundary $\bdry
\widehat X_n$.

\begin{theorem}
$\widehat{X}_n$ equivariantly deformation retracts onto $\widehat{Y}_n$.  Under this retraction, the image $\widehat
Z_n$ of
$\partial\widehat{X}_n$ is the union of $\partial\widehat Y_n$ with the set of graphs in the interior $Y_n$
which have a cut vertex.
\end{theorem}

\begin{proof}  Let $x=(\Gamma_0\supset\ldots\supset\Gamma_k)$ be a point in $\widehat X_n$. 
We define an equivariant deformation retraction $\phi(x,t)$ as follows.  

 If 
$\Gamma_1$ is not the maximal core of $\Gamma_0$, then 
the deformation retraction $\phi$ changes the  metric on $\Gamma_0$  by uniformly 
shrinking all separating edges and rescaling the primary metric on the rest of the graph by a global factor to
retain total volume 1.  The metrics on $\Gamma_k$ for $k\geq 1$ are not affected.  If, on the other hand, $\Gamma_1$ 
is the maximal core  of $\Gamma_0$,
then the deformation retraction shrinks the separating edges in $\Gamma_0$ (i.e. it shrinks $\Gamma_0\setminus \Gamma_1$) 
while simultaneously blowing up the initially degenerate $\Gamma_1$ by an appropriate factor in the primary metric.
In other words, $\Gamma_1$ immediately disappears from the filtration.
  The metrics on
$\Gamma_i$ for
$i\geq 2$ are not affected. 

We now give an explicit formula for $\phi(x,t)$.  Let $m_i$ denote the metric on $\Gamma_i$,  let $S$ denote the
set of separating edges in
$\Gamma_0$, and let $(\Gamma_i)_S$ be the image of $\Gamma_i$ under the map which
collapses each edge in $S$ to a point. The formula depends on whether $m_0(S)=1$ (i.e.
$\Gamma_1$ is the entire core of $\Gamma_0$) or $m_0(S)<1$.

If $m_0(S)<1$, then for $0<t<1$ we have $\phi(x,t)= (\Gamma_0\supset \Gamma_1\supset\ldots\supset \Gamma_k)$, where
the new metric $n_i$ on $\Gamma_i$ is given by
$$n_i(e)=
\begin{cases}
(1+t({{m_0(S)}\over{1-m_0(S)}}))\,m_0(e)&i=0, e\in \Gamma_0\setminus S\\
(1-t)\,m_0(e)& i=0, e\in S\\
m_i(e)& i>0, e\in \Gamma_i\\
\end{cases}$$
For $t=1$, $\phi(x,t)=((\Gamma_0)_S\supset
(\Gamma_1)_S\supset\ldots\supset (\Gamma_k)_S)$, where $(\Gamma_i)_S$ is the image of $\Gamma_i$ under the map which
collapses each edge in $S$ to
a point. The length of each edge $e\in(\Gamma_0)_S$ is $({1\over{1-m_0(S)}})\,m_0(e)$, and the length of
$e\in(\Gamma_i)_S$ is $m_i(e)$.

If $m_0(S)=1$, then for $0<t<1$ we have 
$\phi(x,t)=
(\Gamma_0\supset \Gamma_2\supset\ldots\supset \Gamma_k)$, where the new metric $n_i$ on $\Gamma_i$ is given by
$$n_i(e)=
\begin{cases}
(1-t)m_0(e) + tm_1(e)&i=0,  e\in \Gamma_0\setminus S\\
(1-t)m_0(e) & i=0, e\in S\\
m_i(e)& i>1, e\in \Gamma_i\\
\end{cases}$$
If $t=1$, then $\phi(x,t)=
((\Gamma_1)_S\supset
(\Gamma_2)_S\supset\ldots\supset (\Gamma_k)_S)$, where the 
length of each edge $e$  of $(\Gamma_i)_S$ is equal to  $m_i(e)$ for all $i$.

The deformation retraction restricted to a cell in the case $n=2$ is pictured in Figure~\ref{triangle}. The top line
corresponds to graphs with a degenerate core, and the flow pushes them into strata of $\widehat X_n$  of one
higher dimension. Everywhere else, the flow stays within strata until $t=1$, when the dimension of
the stratum may decrease.

The fact that points in $\bdry\widehat X_n$ land in  $\widehat Z_n$ is
clear.  Now we attack the question of continuity. For this, it will be convenient
to fix a metric on each closed cell of $\hat X_n$. Every top dimensional cell is associated
to a marked trivalent graph, $\Gamma$. Call such a closed cell $\hat{\Sigma}_\Gamma$. 
Let $C$ be a core subgraph of $\Gamma$. For every point of $\hat{\Sigma}_\Gamma$, $C$ has a \emph{level},
which is the unique $i$, such that the metric $m_i$ is defined on $C$ and is not identically zero on $C$.
(If we are looking at a point on $\partial\hat{\Sigma}_\Gamma$ where a subforest has been contracted, then the level
is defined for the image of $C$ under this contraction.)

Let $x,x^\prime$ be two points in $\hat{\Sigma}_\Gamma$, and let $C\subset \Gamma$ be a core subgraph and let $l,
l^\prime$ be the levels
of $C$ in $x$ and $x^\prime$. 
Then define
$$d_C(x,x^\prime) = \sum_{e\in
E(C)}\left|\frac{m_l(e)}{m_l(C)}-\frac{m^\prime_{l^\prime}(e)}{m_{l^\prime}^\prime(C)}\right|.$$ 

Then the metric on on $\hat{\Sigma}_\Gamma$ is defined to be $$d(x,x^\prime)=\sum_{C \subset \Gamma} d_C(x,x^\prime),$$ where the
sum is over all core subgraphs of $\Gamma$ including $\Gamma$ itself. That this metric generates the appropriate topology follows
from Lemma 2.3 of \cite{BF}.

To show that $\phi$ is continuous it suffices to show that on each closed cell $\hat{\Sigma}_\Gamma$, the functions $\phi(x_0,t)$
are
\emph{equicontinuous} as a family indexed by $x_0$ and that each function $\phi(x,t_0)$ is continuous as a function of $x$.
Recall that
\emph{equicontinuous} means that for every $
\epsilon>0 $ there exists a $\delta>0$ such that $|t-s|<\delta \Rightarrow \forall x_0(
d(\phi(x_0,t),\phi(x_0,s))<\epsilon)$.

The continuity of $\phi(x_0,t)$ as a function of $t$ is clear except when $x_0$ represents
$\Gamma\supset\Gamma_1\supset\ldots$, and $\Gamma_1$ is the maximal core of $\Gamma$. 
However, here too $\phi$ is continuous, since by construction of $\widehat{X}_n$, $x_0
:=\underset{t\to 0}{\lim}\,\phi(x_0,t)$.  
 Note that as functions of $t$, the formula for how the length of each edge changes is a
linear map with coefficients bounded by $1$. 
This ensures that the family of functions is equicontinuous, since 
\begin{align*}
d(\phi(x_0,t),\phi(x_0,s)) &=\sum_{C\subset G} d_C(\phi(x_0,t),\phi(x_0,s))\\
&\leq \sum_{C\subset G}\sum_{e\in E(C)}|t-s|\\
&\leq N\cdot |t-s|
\end{align*} where $N$ is a constant independent of $x_0$.

Now we wish to show continuity in $x$. Clearly $\phi$ is  continuous on the interiors of cells, so we need to
consider what happens as we approach the boundary. It will be simpler to analyze what happens as we go from a cell $\Sigma$ to a
codimension one stratum $B$. Suppose $\Sigma$ corresponds to the sequence of graphs $\Gamma_0\supset\Gamma_1\supset\ldots\supset
\Gamma_k$. Then the face $B$ comes from one of two processes. Either it corresponds to contracting an edge $e$:
$(\Gamma_0)_e\supset (\Gamma_1)_e\supset\ldots\supset (\Gamma_k)_e$, or it corresponds to refining the filtration
by inserting a new core graph $C$: $\Gamma_0\supset \ldots\Gamma_i\supset $C$\supset\Gamma_{i+1}\supset \ldots
\supset\Gamma_k$. For every point $x$ on $B$, there is a canonical path, $P_x$, into $\Sigma$. In the first case, it is defined by
expanding the contracted edge in the metric that makes sense, shrinking the other edges in that metric to maintain total volume
$1$. In the second case, the core $C$ is expanded from length zero in the $i$th metric, using a scaled version of the metric that
had been defined on $C$. The other edges in $\Gamma_i$ are scaled down to maintain total volume $1$.

We will show that, for every
$\epsilon>0$ there is a $\delta$ such that $$\forall x\in B\forall y\in P_x (d(x,y)<\delta\Rightarrow
d(\phi(x,t_0),\phi(y,t_0))<\epsilon).$$

This condition will be called \emph{boundary equicontinuity}. This is sufficient to ensure continuity. For example,
to show continuity at a point $z$ on a codimension $2$ face,
let $x$ be a nearby point in the top cell.
Let $y$ be the projection onto one of the nearby codimension $1$ faces (i.e. $x\in P_y$), and $z^\prime$ the projection
onto the codimension 2 face (i.e. $y\in P_{z^\prime}$). Thus if $x$ is sufficiently close to $z$,
then $x,y$ are close, $y,z^\prime$ are close, and $z,z^\prime$ are close. Then
\begin{multline*}
d(\phi(x,t_0),\phi(z,t_0))\leq \\
 d(\phi(x,t_0),\phi(y,t_0))+d(\phi(y,t_0),\phi(z,t_0))
+d(\phi(z,t_0),\phi(z^\prime,t_0)),
\end{multline*}
and by the boundary equicontinuity hypothesis we can make the first two terms uniformly $<\epsilon/3$ and by continuity on
the interior of cells at $z$ we can bound the last term by $\epsilon/3$.

So now let us show boundary equicontinuity. Let $x$ be an interior point and $x^\prime$ be nearby on a codimension $1$ face,
such that $x\in P_{x^\prime}$.

As mentioned above, in one case, $x'=(\Gamma_0)_e\supset\ldots\supset(\Gamma_k)_e$, where the metric on edges  is
unchanged except in the image of the (unique) graph $\Gamma_i$ of the filtration in which $e$ has non-zero length; in
$(\Gamma_i)_S$, edges are scaled by ${1\over{1-m_i(e)}}$. The fact that
$x$ is close to
$x'$ means that
$m_i(e)$ is very small.  

In the other case  $x'=\Gamma_0\supset\ldots\supset \Gamma_i\supset C\supset \Gamma_{i+1}\supset\ldots\supset
\Gamma_k$.  The metric on $C$ is ${1\over m_i(C)}$ times the restriction of $m_i$ to $C$.  The metric on $\Gamma_i$
is $0$ on $C$, and ${1\over {1-m_i(C)}}\,m_i$ on edges not in $C$.  The fact that
$x$ is close to
$x'$  means that
$m_i(C)$ is small.  

It is now routine  to check that $\phi(x,t)$ is uniformly close to $\phi(x',t)$ in all cases.  As an example,
we check one of the more complicated cases, when $x'=\Gamma_0\supset C\supset \Gamma_{1}\supset\ldots\supset
\Gamma_k$, where $C$ is the core of $\Gamma_0$.  Let $|e|$ denote the primary metric on $x$.  Then $|S|+|C|=1$.
The primary length of $e$ in $x'$ is $0$ if $e\in C$ and $|e|/|S|$ if $e\in S$. The secondary length of $e\in C$ in
$x'$ is $|e|/|C|$.  We now compute $\phi(x,t)$ and $\phi(x',t)$ using the formulas above.

Note that $\Gamma_1$ is {\it
not} the core of
$\Gamma_0$, so
$$\phi(x,t)=\big(\Gamma_0\supset \Gamma_1\supset \ldots\supset \Gamma_k\big),$$
for $0<t<1$. The primary length of $e$ is
$$|e|_{\phi(x,t)} = 
\begin{cases}
(1-t) |e|& e\in S\\
(1+t{|S|\over{1-|S|}})|e|=((1-t)|C| + t){|e|\over|C|}&e\in C
\end{cases}$$
On the other hand, 
$C$ {\it is} the core of
$\Gamma_0$ so again $$\phi(x',t)=
\big(\Gamma_0\supset \Gamma_1\supset \ldots\supset \Gamma_k\big),$$
for $0<t<1$.  Now the primary length of $e$ is
$$|e|_{\phi(x',t)} = 
\begin{cases}
(1-t) {|e|\over |S|}& e\in S\\
t{|e|\over|C|}&e\in C
\end{cases}$$

Now, to show that  equicontinuouity at this boundary, we calculate distances.
First, we claim that $d(x,x^\prime)=d_{\Gamma_0}(x,x^\prime)$. 
So let $D\neq \Gamma_0$ be a core subgraph of $\Gamma_0$. Then $D\subset C$. If the primary metric of $x$ vanishes
on $D$, then the first nonvanishing metric is the same for both $x$ and $x^\prime$, and so $d_D(x,x^\prime)=0$.
If $D$ does not vanish in $x$'s primary metric, $D$ vanishes in the primary metric of $x^\prime$, and is rescaled
 in the secondary metric: $m_2^\prime|_D=m_1|_D\cdot \frac{1}{|C|}$. Then 
\begin{align*}
d_D(x,x^\prime)&=
\sum_{e\in E(D)} \left|\frac{m_1(e)}{m_1(D)} - \frac{m^\prime_2(e)}{m^\prime_2(D)}\right| \\
&=\sum_{e\in E(D)}  \left|\frac{m_1(e)}{m_1(D)} - \frac{m_1(e)\cdot |C|^{-1}}{m_1(D)\cdot |C|^{-1}}\right|\\
&=0
\end{align*}

  So we have
\begin{align*}
d(x,x^\prime)&= d_{\Gamma_0}(x,x^\prime)\\
&=\sum_{e\in E(C)}| |e|-|e|^\prime| +\sum_{e\in E(S)}||e|-|e|^\prime| \\
&=
\sum_{e\in E(C)}||e|-0|+\sum_{e\in E(S)} ||e|-\frac{|e|}{|S|}| \\
&=|C|+|S||1-1/|S||\\
&=2|C|.
\end{align*}

On the other hand 
\begin{equation*}
\begin{split}
d(\phi(x,t_0),\phi(x^\prime,t_0))&=d_{\Gamma_0}(\phi(x,t_0),\phi(x^\prime,t_0)) \\
&=\sum_{e\in E(S)}|(1-t_0)|e|-(1-t_0)\frac{|e|}{|S|}|\\
&\phantom{= dummy}+\sum_{e\in E(C)}|(1-t_0)|C|+t)\frac{|e|}{|C|}-t_0\frac{|e|}{|C|}|\\
&=(1-t_0)(|S||1-1/|S||+|C||)\\
&=2(1-t_0)|C|
\end{split}
\end{equation*}

Thus we can take $\delta = \frac{\epsilon}{1-t_0}$, which is independent of $x$.

\end{proof}

\section{Lie bialgebra structure on $\g/\cv$}

Recall that $\mathcal{HG}$ denotes the Hopf algebra spanned by all
oriented graphs (not necessarily connected).
In this section, we will show that the Lie bracket and
cobracket on
$\mathcal{HG}$ introduced in \cite{bialgebra} induce a Lie bracket and
cobracket on $\g/\cv$, and that these are compatible on
Ê $\g/\cv$.

We first recall the definition of the Lie bracket. Let $G$ be a graph,
and let $x$ and $y$ beÊ half-edges of $G$, terminating at the
vertices $v$ and $w$ respectively. Form a new graphÊ as
follows:Ê Cut the edges of $G$ containing $x$ and $y$ in half and glue
$x$ to $y$ to form a new edge $xy$, with vertices $v$ and
$w$. If $y$ was not the other half of $x$ (i.e. $y\neq\bar x$), there
are now two ``dangling" half-edges $\bar x$ and $\bar y$.Ê Glue these
to form another new edge
$\bar x\bar y$.Ê Finally, collapse the edge
$xy$ to a point.Ê We say the resulting graph, denoted $G_{xy}$, is
obtained by {\it contracting the half-edges $x$ and $y$}.

Recall that the Hopf algebra product $G\cdot H$ is the disjoint union
of $G$ and $H$, with appropriate orientation. The bracket
of $G$ and
$H$ is defined to be the sum of all graphs obtained by contracting a
half-edge of $G$ with a half-edge of $H$ in $G\cdot H$:
$$[G,H]=\sum_{x\in G, y\in H} (G\cdot H)_{xy}.$$
For more information about the bracket, we refer to \cite{bialgebra};
there we show, e.g., that there is a second boundary operator on
$\h$, and the bracket measures how far this boundary operator is from
being a derivation.

If $x$ and $y$ belong to separating edges of $G$ and $H$, then $(G\cdot H)_{xy}$
will not be connected, even if $G$ and $H$ are connected.  Thus the bracket on $\mathcal{HG}$
does not restrict to a bracket on $\mathcal{G}.$  It does restrict to a bracket on the subcomplex of 
$\g$ spanned by graphs with no separating edges, but that subcomplex is not quasi-isomorphic to $\g$.
However, we will show that it does
induce a well-defined bracket on $\mathcal{G}/\mathcal{C}$.  The quotient $\mathcal{G}/\mathcal{C}$ has as basis the cosets  $G + \mathcal{C}$,
where $G$ is connected with no cut vertices.
We define the bracket on basis elements
by $[G+\cv,
H+\cv]=[G,H]+\cv$, where $G$ and $H$ are connected with no cut vertices. To see that this is well-defined, we need
the following lemma.

\begin{lemma}  If 
$G,H\neq 0$ are connected and have
no cut vertices, then each term $(G\cdot H)_{xy}$ of $[G,H]$ is connected and has
no cut vertices.  If $G$ or $H$ has a cut vertex, then so does each term
$(G\cdot H)_{xy}$.Ê
\end{lemma}
\begin{proof} 
Since we are restricting to nonzero graphs, we may assume there are no loops at any vertices.
A graph without loops is connected with no cut vertices if and only if there are at least two disjoint
paths between every pair of vertices.
  Let $v$ be the vertex of $G$ adjacent to $x$ and $\bar v$ the vertex adjacent to $\bar x$;
similarly, let $w,\bar w$ be the vertices of
$H$ adjacent to $y, \bar y$.  Choose a path $\alpha$ in $G$ from $v$ to $\bar v$ which does not contain $x$, and $ B $ a path
in $H$ from $w$ to $\bar w$ which does not contain $y$.

If $a$ and $b$ are two vertices of $G$, and one of the two disjoint paths
between them contains $x\bar x$, then we can construct a second
disjoint path in $(G\cdot H)_{xy}$ by replacing $x\bar x$ by $ B $. 
Similarly, if $a$ and $b$ are in $H$, we can replace a path containing
$y\bar y$ by $\alpha$.  If $a$ is in $G$ and $b$ is in $H$, then disjoint
paths can be constructed as follows: To make the first path, join $a$ to
the image of $v$ and the (identical) image of $w$ to $b$; the second path
is obtained by joining $a$ to $\bar v$, then going across 
$\bar x
\bar y$, then joining $\bar w$ to $b$.

If a vertex  is a cut vertex in $G$, its image is a cut vertex in $(G\cdot H)_{xy}$.
\end{proof}

\begin{corollary} The bracket induces a well-defined bracket on $\g/\cv$.
\end{corollary}
\begin{proof} The only subtlety here is that the bracket of two graphs with separating edges (and hence cut vertices) may not be connected.
Let $C_1, C_2\in \cv$, and let $\mathcal{HC}$ be the subspace of $\mathcal{HG}$ spanned by
graphs with cut vertices. Then $[G+C_1, H+C_2]=[G,H]+[G,C_2]+[C_1,H]+[C_1,C_2]\in [G,H]+\mathcal{HC},$  by the lemma. We then appeal to the natural
isomorphism 
$$(\g+\mathcal{HC})/\mathcal{HC}\cong \g/\cv$$
to identify $[G+C_1,H+C_2]$ with $[G,H]$ in $\g/\cv$.
\end{proof}

\begin{remark}  The Lemma shows that the bracket restricts to the subspace of $\g$ spanned by graphs with no cut vertices. 
However, that subspace is not a subcomplex, as the boundary map does not restrict.  This is the reason we are using the quotient
complex $\g/\cv$.
\end{remark}

The cobracket $\mathcal{HG}\to\mathcal{HG}\otimes\mathcal{HG}$ is
defined as follows.Ê We say a pair $\{x,y\}$ of half-edges of a 
graph $G$ is
{\it separating} if the number of components of $G_{xy}$ is greater than that of $G$.Ê
If $G$ is
connected, define
$$\theta(G)=\sum_{\{x,y\} \,\, \text{separating}} A\otimes B + (-1)^a B\otimes
A,$$
where $G_{xy}= A\cdot B$, and $a$ is the number of vertices of $A$.
This gives the coproduct on primitive
elements, and extends to all elements in a standard way (see
\cite{bialgebra}).Ê We have

\begin{lemma} Let $\{x,y\}$ be a separating pair of half-edges in an
oriented graph $G$, with $G_{xy}=A\cdot B$. If $G$ has a cut vertex,
then at least one of $A$ or $B$ has a cut vertex.Ê If $G$ is connected,
then both $A$ and $B$ are connected.
\end{lemma}

\begin{proof} The proof is straightforward.
\end{proof}

Thus the cobracket induces a cobracket on
$(\g+\mathcal{HC})/\mathcal{HC}\cong \g/\cv$ defined on a basis element
$G+\cv$, where
$G$ is a connected graph with no cut vertices, by
$$\theta(G+\cv)=\sum_{\{x,y\} \,separating} (A+\cv)\otimes (B+\cv)+
(-1)^a (B+\cv)\otimes (A+\cv).$$

We now check that the bracket and cobracket are compatible on $\g/\cv$,
making $\g/\cv$ into a Lie bialgebra:

\begin{proposition}Ê The bracket and cobracket satisfy
$$\theta[G+\cv,H+\cv]+[\theta(G+\cv),H+\cv] + (-1)^g
[G+\cv,\theta(H+\cv)]=0,$$
where $g$ is the degree of $G$.
\end{proposition}

\begin{proof} Let $G+\cv$ and $H+\cv$ be basis elements of $\g/\cv$,
i.e. $G$ and $H$ are connected with no cut vertices.
We compute
\begin{align*}
\theta[G+\cv,H+\cv]+[\theta(G+\cv),H+\cv] + (-1)^g
[G+\cv,\theta(H+\cv)]\\
=\theta[G,H]+[\theta(G),H] + (-1)^g
[G,\theta(H)].
\end{align*} Because both $G$ and $H$ have no cut vertices, they also
have no separating edges, so the last sum is zero by Theorem 1
of
\cite{bialgebra}.
\end{proof}

\section{Calculations}

In this section we present our computations of the rational homology of
$\g^{(n)}$ for $n \le 7$, briefly describing the algorithm, but omitting
the raw code. Details for a similar algorithm can be found in \cite{gerlits}.

The program first enumerates all
trivalent graphs with no cut vertices. There is only one such graph with
fundamental group of rank 2, the \emph{theta graph:} two vertices
connected by three edges. If we have a list of all graphs with
fundamental group of rank $n-1$, we can obtain the list for rank $n$ by
applying one of the following two operations to all the graphs in
every possible way. The first operation takes two distinct edges
of the graph, subdivides them by adding a new vertex at the middle
of each, and adds a new edge between the two new vertices. The
second operation adds two new vertices in the interior of a single
edge and connects the two new vertices by a new edge. It is not
hard to see that if $G$ has no cut vertices, then it can be
obtained from a lower rank graph which also has no cut vertices
using one of these two operations.

The same graph will be listed several times. To eliminate the
duplications, we transform each new graph into a \emph{normal
form;} two graphs in normal form are isomorphic if and only if
they are identical. The graph is stored as the adjacency matrix
$a_{ij}$ for $i < j$; i.e., $a_{ij} = k$ if vertex $i$ is
connected to vertex $j$ by $k$ edges. The normal form of the graph
is the ordering of the vertices which yields the matrix latest in
the lexicographic ordering. Permutations of the
vertices are listed, and the matrices are compared. 
The number of permutations needed
is reduced by distinguishing 3 types of vertices:
those contained in a multiple edge, those contained in a triangle,
and the rest. Only vertices of the same type need to be permuted
among themselves. 

Next, we enumerate all graphs of valence 3 or higher, with no cut
vertices and with fundamental group of rank at most 7, by
successively contracting edges of the trivalent graphs. Cut
vertices may develop during this process; in this case the graph
is discarded. Then we examine each graph to see whether it has any
orientation-reversing automorphisms, and if so, discard it. A graph
with an orientation-reversing automorphism is zero in the
graph complex since such a graph is equal to minus itself and
the base field is of characteristic zero.

Finally, we compute the matrix of the boundary map
\begin{equation*}
\partial(G) = \sum_{e} G_e
\end{equation*}
by transforming each contracted graph $G_e$ into normal form and
comparing it to the list of lower rank graphs. The output of the
program is a sparse matrix; its rank was computed by simple
Gaussian elimination in the case of the smaller matrices and by
the software package {\tt scilab} in the case of the largest ones.

Recall that $\g^{(n)}$ denotes the subcomplex of the graph complex spanned by 
rank $n$ graphs, and that $\cv^{(n)}$ is the subcomplex spanned by graphs with cut vertices.
We obtain the following rank $n$ quotient complexes $\mathcal G^{(n)}/\mathcal C^{(n)}$ for values of $n$ less than $8$:

\newcommand{\s}[1]{{\mbox{\small #1}}}

$\g^{(3)}/\cv^{(3)}$:
$$
\begin{CD}
0 @>>> \underset{2}{C_4} @>{1}>> \underset{1}{C_3} @>{0}>> \underset{0}{C_2} @>>> 0
\end{CD}
$$

\vspace{.5em}

$\g^{(4)}/\cv^{(4)}$:
$$
\begin{CD}
0 @>>> \underset{4}{C_6} @>{3}>> \underset{3}{C_5} @>{0}>> \underset{0}{C_4} @>{0}>> \underset{1}{C_3} 
@>{1}>> \underset{1}{C_2} @>>> 0
\end{CD}
$$

\vspace{1em}

$\g^{(5)}/\cv^{(5)}$:
\begin{multline*}
\begin{CD}
0 @>>>\underset{14}{C_8}@>{12}>>\underset{19}{C_7} @>{7}>> \underset{12}{C_6} @>{5}>> \underset{12}{C_5} @>{7}>>
\underset{10}{C_4} @>{3}>>
\underset{3}{C_3}  @>{0}>> 
\end{CD}\\
\begin{CD}
\underset{0}{C_2} @>>> 0
\end{CD}
\end{multline*}

\vspace{1em}

$\g^{(6)}/\cv^{(6)}$:
\begin{multline*}
\begin{CD}
0 @>>>\underset{54}{C_{10}}@>{52}>>\underset{128}{C_9} @>{76}>> \underset{177}{C_8} @>{101}>> \underset{218}{C_7} @>{116}>>
\underset{177}{C_6} @>{61}>>
\underset{72}{C_5}  @>{11}>>
\end{CD}\\
\begin{CD}
\underset{12}{C_4} @>{1}>>
\underset{2}{C_3}  @>{1}>> \underset{1}{C_2} @>>> 0
\end{CD}
\end{multline*}

\vspace{1em}

$\g^{(7)}/\cv^{(7)}$:
\begin{multline*}
\begin{CD}
0 @>>>\underset{298}{C_{12}}@>{295}>>\underset{1123}{C_{11}} @>{828}>> \underset{2388}{C_{10}} @>{1560}>> \underset{3530}{C_9}
@>{1969}>>
\underset{3362}{C_8} @>{1393}>>
\end{CD}\\
\begin{CD}
\underset{1933}{C_7}  @>{540}>> 
\underset{678}{C_6} @>{138}>> \underset{173}{C_5} @>{35}>>
\underset{41}{C_4} @>{6}>>
\underset{6}{C_3}  @>{0}>> \underset{0}{C_2} @>>> 0
\end{CD}
\end{multline*}

The number printed under the chain group $C_i$ is its dimension,
the number printed above the arrow $\partial_i\colon C_i \to
C_{i-1}$ is its rank. Thus we have the following.

\begin{theorem} \label{computation}
The rational homology $H_i(\g^{(n)})$ of the commutative graph
complex is zero for all $2\le n \le 7$ and all $i$ except for
\begin{alignat*}{1}
&H_2(\g^{(2)}) \cong \Q \\
&H_4(\g^{(3)}) \cong \Q \\
&H_6(\g^{(4)}) \cong \Q \\
&H_8(\g^{(5)}) \cong \Q^2 \\
&H_{10}(\g^{(6)}) \cong \Q^2
\qquad H_7(\g^{(6)}) \cong \Q \\
&H_{12}(\g^{(7)}) \cong \Q^{3}
\qquad H_9(\g^{(7)}) \cong \Q.
\end{alignat*}
\end{theorem}

For $2 \le n \le 5$, the computation takes only a few minutes. For
$n = 6$, it took several hours of CPU time, and for $n = 7$,
several thousand hours, even though the elimination of graphs with
cut vertices reduces the size of the computation by about 30\%.

As Kontsevich realized, any metrized Lie algebra produces classes in
trivalent graph homology, and so the abundance of top dimensional
homology is perhaps not surprising. Indeed, these trivalent classes 
correspond to finite type three manifold invariants.
On the other hand, the presence of two codimension $3$ classes is rather 
tantalizing.

The source code for the program is available in the {\tt source Folder} available
with the ``source" for this paper on arXiv.org. Please look at the {\tt readme} file first.
There is also a {\tt data Folder} available in the same place.

\medskip\noindent{\bf Acknowledgements.}
The first author was partially supported by NSF grant
DMS-0305012.
The third author was partially supported by NSF grant DMS-0204185. The computations were done at the NIIFI Supercomputing Centre in Budapest, Hungary.  We thank Craig Jensen for pointing out several typographical errors in
an earlier draft.


\begin{thebibliography}{999}


\bibitem{bgrt} {\bf D. Bar-Natan, S. Garoufalidis, L. Rozansky and
D. Thurston}, \textsl{The Aarhus integral of rational homology 3-spheres
I: A highly nontrivial flat connection on $S^3$.}, Selecta
Math. (N.S.) 8 (2002), no. 3, 315--319.

\bibitem{BF} {\bf Mladen Bestvina and Mark Feighn}, 
\textsl{The topology at infinity of $Out(F_n)$}
Invent. Math. 140 (2000), no. 3, 651--692


\bibitem{fusion} {\bf James Conant},
\textsl{Fusion and fission in graph complexes}, Pac.
J. Math., Vol. 209, No.2 (2003), 219-230.

\bibitem{bialgebra}Ê{\bf James Conant
and Karen Vogtmann},
\textsl{Infinitesimal operations on complexes of graphs}, 
Math. Ann.327 (2003), no. 3, 545--573.

\bibitem{exposition}
ÊÊÊÊ {\bf{James Conant and Karen Vogtmann}},
\textsl{On a theorem of Kontsevich}, Algebr. Geom. Topol. 3 (2003), 1167-1224.

\bibitem{gerlits} {\bf Ferenc Gerlits}, \textsl{Invariants in chain complexes
of graphs}, Cornell Ph.D. thesis, 2002.

\bibitem{gan} {\bf Wee Liang Gan}, \textsl{On a theorem of Conant-Vogtmann}, preprint 2004, math.QA/0404173.

\bibitem{Ko1} {\bf Maxim Kontsevich}, \textsl{Formal (non)commutative
symplectic geometry.} The Gelfand
Mathematical Seminars, 1990--1992, 173--187, Birkh\"auser Boston,
Boston, MA, 1993.

\bibitem{Ko2} {\bf Maxim Kontsevich}, \textsl{Feynman diagrams and
low-dimensional topology.} First European Congress of
Mathematics, Vol. II (Paris, 1992), 97--121, Progr. Math., 120,
Birkh"user, Basel, 1994.

\bibitem{kt} {\bf Greg Kuperberg and Dylan P. Thurston} 
\textsl{Perturbative
3-manifold invariants by cut-and-paste topology}, preprint 1999, UC 
Davis Math 1999-36,  math.GT/9912167

\bibitem{lmo} {\bf T.Q.T. Le, J. Murakami and T. Ohtsuki},
\textsl{On a universal perturbative invariant of $3$-manifolds}, Topology
{\bf 37-3} (1998)

\bibitem{survey} {\bf Karen Vogtmann}, \textsl{Automorphisms of free groups
and Outer Space} Geometriae Dedicata, v. 94: 1--31, 2002.

\end{thebibliography}
\end{document}